\documentclass[12pt]{article}
\usepackage{amssymb}
\usepackage{latexsym}
\usepackage{amsmath}
\def\IND#1{{\mathbb I}_{{\left[ #1 \right]}}}
\def\RR{\mathbb{R}}
\newcommand{\argmin}{\mathop{\arg\!\min}}
\newcommand{\argmax}{\mathop{\arg\!\max}}
\newcommand{\Argmin}{\mathop{\mathrm{Arg}\!\min}}

\newtheorem{theorem}{Theorem}

\newtheorem{definition}{Definition}

\newtheorem{proposition}{Proposition}

\newtheorem{rema}{Remark}
\newenvironment{remark}{\begin{rema} \rm}{\end{rema}}
\newtheorem{exam}{Example}

\newenvironment{proof}[1][Proof]{\textbf{#1.} }{\ \rule{0.5em}{0.5em}}

\begin{document}
\title{Recursive Aggregation of Estimators by Mirror Descent Algorithm with
Averaging}
%
%\titlerunning{Online  Aggregation of Classifiers with a Mirror Descent Algorithm}
%
%
\author{Anatoli Juditsky\\ Laboratoire de Mod\'elisation et Calcul - Universit\'e Grenoble I\\
B.P. 53, 38041 Grenoble, France\\
{\tt anatoli.iouditski@imag.fr} \and Alexander Nazin\\ Institute
of Control Sciences -  Russian Academy of Sciences \\
65, Profsoyuznaya str., GSP-7, Moscow, 117997, Russia \\
{\tt nazine@ipu.rssi.ru} \and
 Alexandre Tsybakov  \and Nicolas Vayatis\\
 Laboratoire de Probabilit\'{e}s et Mod\`{e}les Al\'{e}atoires\\
UMR CNRS 7599 - Universit\'{e}s Paris VI et Paris VII \\
4, place Jussieu, 75252 Paris Cedex, France \\
{\tt \{tsybakov,vayatis\}@ccr.jussieu.fr}}

%
%\authorrunning{Anatoli Juditsky et al.}
%%%% modified list of authors for the TOC (add the affiliations)
%\tocauthor{ Anatoli Juditsky (Universit\'e Grenoble I), Alexander Nazin
%(Institute of Control Sciences of RAS),
% Alexandre Tsybakov (Universit\'e Paris VI), Nicolas Vayatis (Universit\'e Paris VI)}
%
%\institute{Laboratoire de Mod\'elisation et Calcul - Universit\'e Grenoble I\\
%B.P. 53, 38041 Grenoble, France\\
%\email{anatoli.iouditski@imag.fr} \and Institute of
%Control Sciences -  Russian Academy of Science \\
%65, Profsoyuznaya str., GSP-7, Moscow, 117997, Russia \\
%\email{nazine@ipu.rssi.ru} \\
%\and Laboratoire de Probabilit\'es et Mod\`eles Al\'eatoires - Universit\'{e} Paris VI\\
%4, place Jussieu, 75252 Paris Cedex,
%France \\
%\email{\{tsybakov,vayatis\}@ccr.jussieu.fr}}

\maketitle              % typeset the title of the contribution

\begin{abstract}
We consider a recursive algorithm to construct an aggregated
estimator from a finite number of base decision rules in the
classification problem. The estimator approximately mi\-ni\-mi\-zes
a convex risk functional under the $\ell_1^{}$-constraint. It is
defined by a stochastic version of the mirror descent algorithm
(i.e., of the method which performs gradient descent in the dual
space) with an additional averaging. The main result of the paper is
an upper bound
%the rate of convergence
for the expected accuracy of the proposed estimator. This bound is
of the order $\sqrt{(\log M)/t}$ with an explicit and small constant
factor, where $M$ is the dimension of the problem and\, $t$ stands
for the sample size.
%The bound is provided with an explicit and small constant factor.
A similar bound is proved for a more general setting that covers,
in particular, the regression model with squared loss.
\end{abstract}

%%%%%%%%%%%%%%%%%%%%%%%%%%%%%%%%%%%%%%%%%%%%%%%%%%%%%%%%%%%

%\centerline{{\bf\S ~1. Introduction}}\bigskip
\section{Introduction}\label{sec:Introduction}

The methods of Support Vector Machines (SVM) and boosting recently
became widely used in the classification practice (see, {\em e.g.},
\cite{Fre95,Sch90,sfbl,v98}). These methods are based on
minimization of a convex empirical risk functional with a penalty.
Their statistical analysis is given, for instance, in
\cite{BJM03,LuVa04, ScSt03,Zha04} where one can find further
references. In these papers, the classifiers are analyzed as if they
were exact minimizers of the empirical risk functional but in
practice this is not necessarily the case. Moreover, it is assumed
that the whole data sample is available, but often it is interesting
to consider the online setting where the observations come
one-by-one, and recursive methods need to be implemented.

There exists an extensive literature on recursive classification
starting from Perceptron and its various modifications (see, {\em
e.g.}, the monographs \cite{abr70a,abr70b,Tsyp70} and the related
references therein, as well as the overviews in \cite{CCG04,dgl}).
We mention here only the methods which use the same loss functions
as boosting and SVM, and which may thus be viewed as their online
analogues. Probably, the first technique of such kind is the method
of potential functions, some versions of which can be considered as
online analogues of SVM (see \cite{abr70a,abr70b} and \cite{dgl},
Chapter~10). Recently, online analogues of SVM and boosting-type
methods using convex losses have been proposed in \cite{KSW04,
Zha04b}. In particular, in \cite{Zha04b}, a stochastic gradient
algorithm with averaging is studied for a general class of loss
functions (cf. \cite{pj92}). All these papers use the standard
stochastic gradient method for which the descent takes place in the
initial parameter space.

In this paper, we also suggest online versions of boosting and SVM,
but based on a different principle: the gradient descent is
performed in the dual space. Algorithms of this kind are known as
mirror descent methods \cite{NeYu83}, and they were initially
introduced for deterministic optimization problems. Their advantage,
as compared to the standard gradient methods, is that the
convergence rate depends logarithmically on the dimension of the
problem. Therefore, they turn out to be very efficient in
high-dimensional problems \cite{BMN01}.

Some versions of the original mirror descent method of Nemirovski
and Yudin \cite{NeYu83}  were derived independently in the learning
community and have been applied to classification and other learning
problems in the papers \cite{HeKiWa99,KiWa97,KiWa01} where bounds
for the relative risk criterion were obtained. However, these
results are formulated in a deterministic setting and they do not
extend straightforwardly to the standard stochastic analysis with a
mean risk criterion (see \cite{KiWa97,CBG04,CCG04b} for insights on
the connections between the two types of results). Below we propose
a novel version of the mirror descent method which attains the
optimal bounds of the mean risk accuracy. Its main difference from
the previous methods is the additional step of averaging of the
updates.

The goal of this paper is to construct an aggregated decision rule:
we introduce a fixed and finite base class of decision functions,
and we choose the weights in their convex or linear combination in
an optimal way. The optimality of weights is understood in the sense
of minimization of a convex risk function under the
$\ell_{1}^{}$-constraints on the weights. This aggregation problem
is similar to those considered, for instance, in \cite{JuNe00} and
\cite{Tsy03} for the regression model with squared loss. To solve
the problem, we propose a recursive algorithm of mirror descent type
with averaging of the updates. We prove that the algorithm converges
with a rate of the order $\sqrt{(\log M)/t}$, where $M$ is the
dimension of the problem, and $t$ stands for the sample size.
%Moreover, we provide a bound with an explicit and small constant factor.

The paper is organized as follows. First, we give the problem
statement and formulate the main result on the convergence rate
(Section \ref{sec:ProblemResult}). Then, the algorithm is described
(Section \ref{sec:mirror}) and the proof of the main result is given
(Section \ref{sec:rate}). In Section \ref{sec:gener}, the result is
extended to general loss functions and to general estimation
problems. Discussion is given in Section \ref{sec:conclu}.

\bigskip
%%%%%%%%%%%%%%%%%%%%%%%%%%%%%%%%%%%%%%%%%%%%%%%%%%%%%%%%%%

%\centerline{{\bf\S ~2. }}\bigskip
\section{Set-up and Main Result}\label{sec:ProblemResult}

We consider the problem of binary classification. Let $(X,Y)$ be a
pair of random variables with values in $\mathcal{X}\times\{-1,+1\}$
where $\mathcal{X}$ is a feature space. A decision rule
$g_{f}^{}:\mathcal{X}\to\{-1,+1\}$, corresponding to a measurable
function $f:\mathcal{X}\to{\mathbb{R}}$ is defined as $g_{f}^{}(x) =
2\IND{f(x)>0}-1$, where $\IND{\cdot}$ denotes the indicator
function. A standard measure of quality of a decision rule
$g_{f}^{}$ is its risk which equals to the probability of
misclassification $R(g_{f}^{}) =\mathbb{P}\{Y\neq g_{f}^{}(X)\}$.
Optimal decision rule is defined as $g_{f^*}^{}$\,, where $f^*$ is a
minimizer of $R(g_f^{})$ over all measurable $f$. The optimal rule
is not implementable in practice since the distribution of $(X,Y)$
is unknown. In order to approximate $g_{f^*}^{}$, one looks for
empirical decision rules $\widehat{g}_{n}$ based on a sample $(X_1,
Y_1), \dots , (X_{n}, Y_{n})$, where $(X_i,Y_i)$ are independent
random pairs having the same distribution as $(X,Y)$.

An abstract approach to construction of empirical decision rules
\cite{dgl,v98,vc} prescribes to search $\widehat{g}_{n}$ in the form
$\widehat{g}_n=g_{\widehat{f}_n}^{}$\,, where $\widehat{f}_n$ is a
minimizer of the empirical risk (empirical classification error):
\begin{equation}\label{Vst1}
    R_n(g_f^{}) =\frac{1}{n} \sum_{i=1}^n \IND{Y_i\neq g_f^{}(X_i)}
\end{equation}
over all $f$ from a given class of decision rules. Conditions of
statistical optimality of the method of minimization of empirical
classification error (\ref{Vst1}) have been extensively studied (see
in particular \cite{vc, dgl, v98}). However, this method is not
computationally tractable, since the risk functional $R_n$ in
(\ref{Vst1}) is not convex or even continuous. In practice,
efficient methods like SVM and boosting implement numerical
minimization of convex empirical risk functionals different from
(\ref{Vst1}) as it has been first noticed in \cite{Bre98} and
\cite{FHT00}. Theoretical analysis was provided recently in several
papers \cite{BJM03,LuVa04,Zha04, ScSt03} where consistency and rates
of convergence of convex risk minimization methods are established
in terms of the probability of misclassification.

A key argument used in these works is that, under rather general
assumptions, the optimal decision rule $g_{f^*}^{}$ coincides with
$g_{f^A}$ where $f^A$ is optimal decision function in the sense that
it minimizes a convex risk functional called the $\varphi$-risk and
defined by
\[
A (f) = \mathbb{E}\, \varphi(Y f(X))
\]
where $\varphi \, : \mathbb{R} \to \mathbb{R}_+$ is a convex loss
function and $\mathbb{E}$ denotes the expectation. Typical choices
of loss functions are the hinge loss function
$\varphi(x)=(1-x)_+^{}$ (used in SVM), as well as the exponential
and logit losses $\varphi(x)=\exp(-x)$ and
$\varphi(x)=\log_2(1+\exp(-x))$ respectively (used in boosting).

Thus, to find an empirical decision rule $\widehat{g}_{n}$ which
approximates the optimal $g_{f^*}$, we can consider minimizing the
empirical $\varphi$-risk
\[
    A_n(f) =\frac{1}{n} \sum_{i=1}^n \varphi(Y_i f(X_i))~,
\]
which is an unbiased estimate for $A(f)$. This strategy is further
justified by a result in \cite{Zha04} generalized in \cite{BJM03}.
This minimization problem is simpler than the original one because
it can be solved by standard numerical procedure, the functional
$A_n$ being convex. When relevant penalty functions are added, it
leads to some versions of boosting and SVM algorithms. At the same
time, one needs the whole sample $(X_1,Y_1),\dots,(X_n,Y_n)$ for
their implementation, i.e., these are \emph{batch procedures}.

In this paper, we consider the problem of minimization of the
$\varphi$-risk $A$ on a parametric class of functions $f$ when the
data $(X_i,Y_i)$ come sequentially (online setting).

Let us introduce the parametric class of functions in which $f$ is
selected. Suppose that a finite set of base functions $\{h_1,
\dots,h_M\}$ is given, where $h_j \, : \, \mathcal{X} \to [-K, K], \
j=1, \dots, M$, $K>0$ is a constant, and $M\ge 2$. We denote by $H$
the vector function whose components are the base functions:
\begin{equation}\label{eq:base}
\forall x\in\mathcal{X}, \quad    H(x)=\left(h_1(x),
\dots,h_M(x)\right)^T ~.
\end{equation}
A typical example is the one where the functions $h_j$ are decision
rules, i.e., they take values in $\{-1,1\}$. Furthermore, for a
fixed $\lambda>0$ we denote the $\lambda$-simplex in $\mathbb{R}^M$
by $\Theta_{M,\,\lambda}$\,:
\[
\Theta_{M,\,\lambda} = \left\{ \theta = (\theta^{(1)}, \dots,
\theta^{(M)})^T \in \mathbb{R}^M_+\,:\; \sum_{i=1}^M \theta^{(i)}
=\lambda
        \right\} ~.
\]
Introduce a family of $\lambda$-convex combinations of functions
$h_1, \dots,h_M$\,:
\[
\mathcal{F}_{M,\,\lambda}= \left\{f_\theta =\theta^T H  \, : \,
\theta\in\Theta_{M,\,\lambda} \right\} ~.
\]
over which minimization of the $\varphi$-risk $A$ will be performed.
The minimization of $A(f)$ over all $f\in \mathcal{F}_{M,\,\lambda}$
is equivalent to the minimization of $A(f_{\theta})$ over all
$\theta \in \Theta_{M,\,\lambda}$\,, so we simplify the notation and
write in what follows
\[
A(\theta)\triangleq A(f_\theta)~.
\]
Define the vector of optimal weights of the $\lambda$-convex combination of the base
functions as a solution to the minimization problem
\begin{equation}\label{A}
\min_{\theta\in\Theta_{M,\,\lambda}} A(\theta).
\end{equation}
We assume that the distribution of $(X,Y)$ is unknown, hence the
function $A$ is also unknown, and its direct minimization is
impossible. However, we have access to a training sample of
independent pairs $(X_i, Y_i)$, having the same distribution as
$(X,Y)$ that are delivered sequentially and may be used for
estimation of the optimal weights.

In the following section, we propose a stochastic algorithm based on
the mirror descent principle which, at the $t$-th iteration, yields
the estimate $\widehat{\theta}_t= \widehat{\theta}_t((X_1, Y_1),
\dots , (X_{t-1}, Y_{t-1}))$ of the solution to the problem
(\ref{A}). The estimate $\widehat{\theta}_t$ is measurable with
respect to $(\widehat{\theta}_{t-1},X_{t-1}, Y_{t-1})$, which means
that the algorithm fits with the online setting. In order to obtain
the updates of the algorithm, it is sufficient to have random
realizations of the sub-gradient of $A$ which have the form:
\begin{equation}\label{1}
u_{i} (\theta) = \varphi\,'(Y_i \theta^T H(X_i)) Y_i H(X_i) \in
\RR^M, \quad  i=1,2,\dots,
\end{equation}
where $\varphi\,'$ represents an arbitrary monotone version of the
derivative of $\varphi$ (one may take, for instance, the right
continuous version).

Given $\widehat{\theta}_t$, the convex combination
$\widehat{\theta}_t^{\;T}H(\cdot)$ of the base functions can be
constructed,  and it defines an aggregated decision rule
\[
\tilde g_t(x) = 2\IND{\widehat{\theta}_t^{\;T}H(x)>0}-1~.
\]
Statistical properties of this decision rule are described by the
following result which establishes the convergence rate for the
expected accuracy of the estimator $\widehat{\theta}_t$ with respect
to the $\varphi$-risk.

\begin{theorem}\label{thm}
For a given convex loss function $\varphi$, for a fixed number
$M\ge 2$ of %classifiers
base functions
and a fixed value of $\lambda > 0$, let
the estimate
$\widehat{\theta}_t$ be defined by the algorithm of Subsection
\ref{sub:tuning}. Then, for any integer $t\ge 1$,
\begin{equation}\label{eq:th1}
\mathbb{E}\, A (\widehat{\theta}_t) - \min_{\theta\in\Theta_{M,
\lambda}} A(\theta) \le C \,
%\sqrt{\frac{\ln M}{t}}
\frac{(\ln M)^{1/2}\sqrt{t+1}}{t} \,,
\end{equation}
where $C = C(\varphi, \lambda) = 2\, \lambda L_{\varphi}(\lambda)$
and $L_{\varphi}(\lambda) = K \sup_{x \in [-K\lambda, K\lambda]}
|\varphi'(x)|$.
\end{theorem}

For example, Theorem \ref{thm} holds with constant $C=2\,$
in a typical case where we deal with convex ($\lambda=1$)
aggregation of base classifiers $h_j$ taking values in $\{-1,1\}$
and we use the hinge loss $\varphi(x)=(1-x)_+$. We also note that
Theorem~\ref{thm} is distribution free: there is no assumption on
the joint distribution of $X$ and $Y$ except, of course, that $Y$
takes values in $\{-1,1\}$ since we deal with the classification
problem.

%\begin{remark}({\sc the algorithm does not implement the ERM.})
%Note that in this setting, there is no need to introduce the
%empirical counterpart of the target risk functional. The estimator
%is NOT an exact or approximate (recursively obtained) empirical risk
%minimizer, but rather an approximate true risk minimizer based on
%the observations of the stochastic gradient.
%\end{remark}

\begin{remark}({\sc efficiency.}) The rate of convergence of order $\sqrt{(\ln M)/t}$
is typical without low noise assumptions (as they are introduced in
\cite{Tsy04}). Batch procedures based on minimization of the
empirical convex risk functional present a similar rate. From the
statistical point of view, there is no remarkable difference between
batch and our mirror descent procedure. On the other hand, from the
computational point of view, our procedure is quite comparable with
the direct stochastic gradient descent. However, the mirror descent
algorithm presents two major advantages as compared both to batch
and to direct stochastic gradient: (i) its behavior with respect to
the cardinality of the base class is better than for direct
stochastic gradient descent (of the order of $\sqrt{\ln M}$ in
Theorem \ref{thm}, instead of $M$ or $\sqrt{M}$ for direct
stochastic gradient); (ii) mirror descent presents a higher
efficiency especially in high-dimensional problems since its
algorithmic complexity and memory requirements are of strictly
smaller order than for corresponding batch procedures (see
\cite{JuNe00} for a comparison).
\end{remark}

\begin{remark}({\sc optimality of the convergence rate.}) Using
the techniques of \cite{JuNe00} and \cite{Tsy03} it is not hard to
prove the minimax lower bound on the excess risk $\mathbb{E}\, A
(\widehat{\theta}_t) - \min_{\theta\in\Theta_{M, \lambda}}
A(\theta)$ having the order $\sqrt{(\ln M)/t}$ for $M\ge
t^{1/2+\delta}$ with some $\delta>0$. This indicates that the
upper bound of Theorem~\ref{thm} is rate optimal for such values
of $M$.
\end{remark}

\begin{remark}({\sc choice of the base class.})
We point out that the good behavior of this method crucially relies
on the choice of the base class of functions $\{h_j\}_{1 \le  j \le
M}$. A natural choice would be to consider a symmetric class in the
sense that if an element $h$ is in the class, then $-h$ is also in
the class. For a practical implementation, some initial data set
should be available in order to pre-select a set of $M$ functions
(or classifiers) $h_j$. Another choice which is practical and widely
spread is to choose very simple base elements $h_j$ such as decision
stumps; nevertheless, aggregation can lead to good performance if
their cardinality $M$ is very large. As far as theory is concerned,
in order to provide a complete statistical analysis, one should
establish approximation error bounds on the quantity $\inf_{f \in
\mathcal{F}_{M, \lambda}} A(f) - \inf_f A(f)$ showing that
the richness of the base class is reflected both by diversity
(orthogonality or independence) of the $h_j$'s and by its
cardinality $M$. For example, one can take $h_j$'s as the
eigenfunctions associated to some positive definite kernel. We refer
to \cite{TaVa04} for related results, see also \cite{ScSt03}. The
choice of $\lambda$ can be motivated by similar considerations. In
fact, if the approximation error is to be taken into account, it
might be useful to take $\lambda$ depending on the sample size $t$
and tending to infinity with some slow rate ({\em cf.}
\cite{LuVa04}). A balance between the stochastic error as given in
Theorem \ref{thm} and the approximation error would then determine
the optimal choice of $\lambda$.  These considerations are left
beyond the scope of the paper, since we focus here on the
aggregation problem.
\end{remark}

\bigskip

\section{Definition and Discussion of the Algorithm}\label{sec:mirror}

In this section, we introduce the proposed algorithm. It is based
on the mirror descent idea going back to Nemirovski and Yudin
\cite{NeYu83} and is the stochastic counterpart of Nesterov's  primal-dual
subgradient method of deterministic convex optimization, studied in \cite{Nest05-1} and \cite{Nest05-2}.  We first give some definitions and recall some
facts from convex analysis.
\medskip

\subsection{Proxy functions}\label{sub:proxy}
We will denote by $E=\ell^M_1$ the space $\mathbb{R}^M$ equipped
with the 1-norm
$$
\|z\|_1^{}=\sum_{j=1}^M |z^{(j)}|
$$
and by $E^* =\ell^M_{\infty}$ the dual space which is $\mathbb{R}^M$
equipped with the sup-norm
\[
 \| z \|_\infty^{} = \max_{\|\theta\|_1^{} =1}
  z^T\theta = \max_{1 \le j \le M} |z^{(j)}| \, , \quad\forall\, z\in E^*,
\]
with the notation $z=(z^{(1)}, \dots, z^{(M)})^T$.

Let $\Theta$ be a convex, closed set in $E$. For a given
 parameter $\beta>0$ and a convex function $V : \Theta \to
\mathbb{R}$, we call {\em $\beta$-conjugate} function of $V$ the
Legendre-Fenchel type transform of $\beta V$:
\begin{equation}\label{fench}
\forall\, z\in E^* , \quad W_\beta(z) = \sup_{\theta\in\Theta} \left\{ -z^T
\theta -\beta V(\theta)\right\} \, .
\end{equation}

Now we introduce the key assumption (Lipschitz condition in dual
norms $\|\cdot\|_1^{}$ and $\|\cdot\|_\infty^{}$) that will be
used in the proofs of Theorems \ref{thm} and \ref{th2}.

\bigskip

\noindent {\bf Assumption (L).} \textit{A convex function $V: \Theta
\to \mathbb{R}$ is such that its $\beta$-conjugate $W_\beta$ is
continuously differentiable on $E^*$ and its gradient $\nabla
W_{\beta}$ satisfies:
\[
\| \nabla W_{\beta}(z) - \nabla W_\beta(\,\tilde{z}\,) \|_1^{} \le
\frac{1}{\alpha\beta} \|z-\tilde{z} \|_\infty^{} \,, \quad \forall\, z,
\tilde{z} \in E^*,  \;\beta>0,
\]
where $\alpha>0$ is a constant independent of $\beta$.}

\medskip

This assumption is related to the notion of strong convexity w.r.t.
the $\|\cdot \|_1$-norm (see, e.g., \cite{BeNe99,RoWe98}).

\begin{definition}
Fix $\alpha>0$. A convex function $V: \Theta \to \mathbb{R}$
 is said to be {\em $\alpha$-strongly convex}  with respect to the norm
 $\|\cdot \|_1$ if
\begin{equation}\label{sconv}
V(sx+(1-s)y) \le sV(x)+(1-s)V(y) - \frac{\alpha}{2} s(1-s) \|x-y\|_1^2
\end{equation}
for all $x, y \in \Theta$ and any $s\in [0,1]$.
\end{definition}

The following proposition sums up some properties of
$\beta$-conjugates and, in particular, yields a sufficient condition
for Assumption (L).

\begin{proposition}\label{prop:conj}
Consider a convex function $V : \Theta \to \mathbb{R}$ and a
strictly positive parameter $\beta$. Then, the $\beta$-conjugate
$W_{\beta}$ of $V$ %satisfies
has
the following properties.
\begin{enumerate}
\item The function $W_\beta : E^*\to \mathbb{R}$ is convex and has a
conjugate $\beta V$, i.e.
\[
\forall\, \theta \in \Theta , \quad \beta V(\theta) = \sup_{z\in E^*} \left\{
-z^T \theta - W_\beta(z) \right\} \, .
\]
\item If~ $V$ is $\alpha$-strongly convex with respect to the norm $\| \cdot \|_1$
 then
\begin{itemize}
\item[(i)] \  Assumption (L) holds true,
\item[(ii)] $ \displaystyle \
\argmax_{\theta\in\Theta} \left\{-z^T \theta -\beta V(\theta)
        \right\} = -\nabla W_\beta(z) \in\Theta~.
$
\end{itemize}
\end{enumerate}
\end{proposition}

For a proof of this proposition we refer to \cite{BeNe99,RoWe98}.
Some elements of the proof are given in the Appendix, subsection B.

\medskip

\begin{definition}\label{proxyf}
We call a function $V : \Theta \rightarrow \RR_+$ {\em
proxy function} if it is convex, and
\begin{itemize}
\item[(i)]\ there exists a point $\theta_*\in\Theta$ such that $ \min_{\theta\in
\Theta} V(\theta)=V(\theta_*)~, $
\item[(ii)] \  Assumption (L) holds true.
\end{itemize}
\end{definition}

\medskip

\noindent {\bf Example.} Consider the entropy type proxy function:
\begin{equation}\label{2}
\forall\, \theta \in \Theta_{M,\,\lambda}\, , \quad V(\theta) = \lambda \ln
\left(\frac{M}{\lambda}\right) + \sum_{j=1}^M \theta^{(j)} \ln \theta^{(j)} \,
\end{equation}
(where $0\ln 0 \triangleq 0$) which has a single minimizer $\theta_*
= (\lambda/M, \dots, \lambda/M)^T$ with $V(\theta_*)=0$. It is easy
to check (see Appendix, subsection B) that this function is $\alpha$-strongly
convex with respect to the norm $\| \cdot \|_1$, with the parameter
$\alpha=1/\lambda$. An important property of this choice of $V$ is
that the optimization problem (\ref{fench}) can be solved explicitly
so that $W_\beta$ and $\nabla W_\beta$ are given by the following
formulas:
\begin{eqnarray}\label{Wbeta}
  \forall\, z\in E^* , \quad W_{\beta} (z) &=& \lambda
\beta\,\ln\left(\frac{1}{M}\sum_{k=1}^M
     e^{-z^{(k)}/\beta} \right),
\\ \label{gibbs}
 \frac{\partial W_\beta(z)}{\partial z^{(j)}}
\displaystyle &=& - \lambda e^{-z^{(j)}/\beta} \left(\displaystyle\sum_{k=1}^M
e^{-z^{(k)}/\beta}\right)^{-1}, \; j=1, \dots, M.
\end{eqnarray}
Assumption (L) for function  (\ref{2})
holds true, as can be easily proved by direct
calculations
without resorting to Proposition \ref{prop:conj} (see Appendix,
subsection A). Furthermore, note that for $\lambda=1$ the following
holds true:
\begin{itemize}
\item the entropy type proxy function as defined in (\ref{2}) corresponds
to the Kullback information divergence between the uniform
distribution on the set $\{1,\dots,M\}$ and the distribution on
the same set defined by probabilities $\theta^{(j)},\, j=1,\dots,M
$,
\item in view of (\ref{gibbs}), the components of the vector $-\nabla
W_\beta(z)$ define a Gibbs distribution on the coordinates of $z$,
with $\beta$ being interpreted as a temperature parameter.
\end{itemize}

\medskip

\subsection{Algorithm}\label{sub:algo}

Mirror descent algorithms are optimization procedures achieving a
stochastic gradient descent in the dual space. The proposed
algorithm presents two modifications: first, it uses updates of the
{\em stochastic} sub-gradient, and also it presents an averaging
step of the iterate outputs. At each iteration $i$, a new data point
$(X_i, Y_i)$ is observed and there are two updates:
\begin{itemize}
\item one is the variable $\zeta_i$ which is defined by the
stochastic sub-gradients $u_{k}(\theta_{k-1})$, $k=1,\dots,i$, as
the result of the descent in the dual space $E^{*}$,
\item the other update is the parameter $\theta_i$ which is the ``mirror
image'' of $\zeta_i$ in the initial space $E$.
\end{itemize}

In order to tune the algorithm properly, we will also need two fixed
positive sequences $(\gamma_i)_{i \ge 1}$ (step size) and
$(\beta_i)_{i \ge 1}$ (``temperature'') such that
$\beta_i\ge\beta_{i-1}$\,, $\forall i\geq 1$. The algorithm is
defined as follows:

\medskip

\bigskip

\fbox{
\begin{minipage}{0.9\textwidth}

\begin{itemize}
\item[$\bullet$] Fix the initial values $\theta_0 \in \Theta$ and $\zeta_0 =0 \in
\RR^M$.
\item[$\bullet$] For $i=1, \dots, t-1$, do the recursive update
\begin{eqnarray}\label{alg1}
\begin{array}{rcl}
\zeta_i &=& \zeta_{i-1} +\gamma_i u_i(\theta_{i-1}) \,,
\\
&& \\
\theta_i &=& -\nabla W_{\beta_i}(\zeta_i)~.
\end{array}
\end{eqnarray}
\item[$\bullet$]  Output at iteration $t$ the following convex combination:
\begin{eqnarray}\label{alg2}
\widehat{\theta}_t = \frac{ \sum_{i=1}^t \gamma_i\theta_{i-1}}{\sum_{i=1}^t
\gamma_i}\,.
\end{eqnarray}
\end{itemize}

\end{minipage}
}

\bigskip

\medskip

Note that if $V$ is the entropy type proxy function defined in
(\ref{2}), the components $\theta^{(j)}_i$ of vector $\theta_i$ from
(\ref{alg1}) have the form
\begin{eqnarray}
\theta^{(j)}_i \displaystyle &=& \frac{\lambda \displaystyle\exp{\left(-
\beta_i^{-1}\sum_{m=1}^i\gamma_m
u_{m,\,j}(\theta_{m-1})\right)}}{\displaystyle\sum_{k=1}^M
\exp{\left(- \beta_i^{-1}\sum_{m=1}^i\gamma_m
u_{m,\,k}(\theta_{m-1})\right)}} \, , \label{EMDA}
\end{eqnarray}
where $u_{m,\,j}(\theta)$ represents the $j$-th component of vector
$u_{m}(\theta)$, $j=1,\dots,M$.
\medskip

\subsection{Heuristic considerations}\label{sub:heur}

Suppose that we want to minimize a convex function $\theta \mapsto
A(\theta)$ over a convex set $\Theta$. If $\theta_0,\dots,\theta_{t-1}$
are the available search points at iteration $t$,
we can provide the affine approximations $\phi_i$ of the
function $A$ defined, for $\theta\in\Theta$, by
\[
\phi_{i}(\theta)= A(\theta_{i-1})+ (\theta -\theta_{i-1})^T \nabla
A(\theta_{i-1}), \;\;\;\; i=1,\dots,t\,.
\]
Here $\theta \mapsto\nabla A(\theta)$ is a vector function belonging
to the sub-gradient of $A(\cdot)$. Taking a convex combination of
the $\phi_i$'s, we obtain an averaged approximation of $A(\theta)$:
\[
\bar\phi_t(\theta)=\frac{\sum_{i=1}^{t}\gamma_i
\left(A(\theta_{i-1})+ (\theta -\theta_{i-1})^T \nabla
A(\theta_{i-1})\right)}{\sum_{i=1}^{t} \gamma_i}.
\]
At first glance, it would seem reasonable to choose as the next
search point a vector $\theta\in \Theta$ minimizing the
approximation $\bar{\phi}_t$, i.e.,
\begin{eqnarray}\label{ai3}
\theta_t=\argmin_{\theta\in \Theta}\bar{\phi}_t(\theta)=\argmin_{\theta\in
\Theta} \, \theta^T\left(\sum_{i=1}^{t}\gamma_i \nabla A(\theta_{i-1})\right).
\end{eqnarray}
However, this does not make any progress, because our
approximation is ``good'' only in the vicinity of search points
$\theta_0,\dots,\theta_{t-1}$. Therefore, it is necessary to
modify the criterion, for instance, by adding a special penalty
$B_t(\theta,\theta_{t-1})$ to the target function in order to keep
the next search point $\theta_t$ in the desired region. Thus, one
chooses the point:
\begin{equation}\label{heur}
\theta_t=\argmin_{\theta\in \Theta}\left[ \theta^T\left(\sum_{i=1}^{t}\gamma_i
\nabla A(\theta_{i-1})\right)+B_t(\theta,\theta_{t-1})\right].
\end{equation}
Our algorithm corresponds to a specific type of penalty
$B_t(\theta,\theta_{t-1}) = \beta_tV(\theta)$, where $V$
is the proxy function.

Also note that in our problem the vector-function $\nabla A(\cdot)$
is not available. Therefore, we replace in (\ref{heur}) the unknown
gradients $\nabla A(\theta_{i-1})$ by the observed stochastic
sub-gradients $u_i(\theta_{i-1})$. This yields a new definition of
the $t$-th search point:
\begin{equation}
\theta_t = \argmin_{\theta\in \Theta}\left[
\theta^T\left(\sum_{i=1}^{t}\gamma_i u_i(\theta_{i-1})\right)+\beta_t
V(\theta)\right] =\argmax_{\theta\in \Theta}\left[ -\zeta_t^T\theta -\beta_t
V(\theta)\right], \label{heur1}
\end{equation}
where
$$
\zeta_t=\sum_{i=1}^{t}\gamma_i u_i(\theta_{i-1}).
$$
Observe that by Proposition \ref{prop:conj}, the solution to this
problem reads as $-\nabla W_{\beta_{t}}(\zeta_t)$ and it is now
easy to deduce the iterative scheme (\ref{alg1}).

\medskip

\subsection{A particular instance of the algorithm}\label{sub:tuning}

We now define the special case of the mirror descent method with
averaging for which Theorem \ref{thm} is proved. We consider the
algorithm described in Subsection \ref{sub:algo} with the entropy
type proxy function $V$ as defined in (\ref{2}) and with the
following specific choice of the sequences $(\gamma_i)_{i \ge 1}$
and $(\beta_i)_{i \ge 1}$:
\begin{equation}\label{gb1}
\gamma_i\equiv 1\,,\quad \beta_i =\beta_0\sqrt{i+1}\,,\quad i=1,2,\dots,
\end{equation}
where
\begin{equation}\label{gb2}
\beta_0 =  \, L_{\varphi}(\lambda)(\ln M)^{-1/2}~.
\end{equation}
Thus, the algorithm becomes simpler and can be implemented in the
following recursive form:
\begin{eqnarray}
 \zeta_i &=& \zeta_{i-1} + u_i(\theta_{i-1})~,
 \label{alg34a} \\
\theta_i &=& -\nabla W_{\beta_i}(\zeta_i)~,
 \label{alg34b}  \\
\widehat{\theta}_i &=& \widehat{\theta}_{i-1}
-\frac{1}{i}\left(\widehat{\theta}_{i-1} - \theta_{i-1}\right),
\quad i=1,2,\dots,
 \label{alg34c}
\end{eqnarray}
with initial values $\zeta_0=0$\,, $\theta_0\in\Theta$\,
and
$(\beta_i)_{i \ge 1}$ from (\ref{gb1}), (\ref{gb2}).

\medskip

\subsection{Comparison to other Mirror Descent Methods}\label{sub:compar}

The versions of mirror descent method proposed in \cite{NeYu83} are
somewhat different from our iterative scheme (\ref{alg1}). One of
them, which is the closest to (\ref{alg1}), is studied in detail in
\cite{BeNe99}. It is based on the recursive relation
\begin{eqnarray}\label{MDA}
\theta_i =-\nabla W_{1}\Big(-\nabla V(\theta_{i-1})+\gamma_i
u_i(\theta_{i-1})\Big), \quad i=1,2,\dots,
\end{eqnarray}
where the function $V$ is strongly convex with respect to the norm
of initial space $E$ (which is not necessarily the space
$\ell^M_1$) and $W_1$ is the 1-conjugate function to $V$.

If $\Theta=\RR^M$ and $V(\theta)=\frac{1}{2}\|\theta\|_2^2$, the scheme of
(\ref{MDA}) coincides with the ordinary gradient method.

For the unit simplex $\Theta=\Theta_{M,1}$ and the entropy type
proxy function $V$ from (\ref{2}), the components $\theta^{(j)}_i$
of vector $\theta_i$ from (\ref{MDA}) are:
\begin{eqnarray}
\nonumber \forall j=1,\dots,M , \quad   \theta^{(j)}_i \displaystyle
&=& \frac{ \theta^{(j)}_{i-1} \exp{(-\gamma_i
u_{i,\,j}(\theta_{i-1}))}}{\displaystyle \sum_{k=1}^M
\theta^{(k)}_{i-1} \exp{(-\gamma_i u_{i,\,k}(\theta_{i-1}))}}
 \\
&=& \frac{\displaystyle \theta^{(j)}_{0} \exp{\left(-
\sum_{m=1}^i\gamma_m u_{m,\,j}(\theta_{m-1})\right)}}{
\displaystyle\sum_{k=1}^M \theta^{(k)}_{0}\exp{\left(-
\sum_{m=1}^i\gamma_m u_{m,\,k}(\theta_{m-1})\right)}}~. \label{EGA}
\end{eqnarray}
The algorithm (\ref{EGA}) is also known as the exponentiated
gradient (EG) method  \cite{KiWa97}. The differences between the
algorithms (\ref{MDA}), (\ref{EGA}) and ours are the following:
\begin{itemize}
\item the initial iterative scheme (\ref{alg1}) is different
from that of (\ref{MDA}), (\ref{EGA}), in particular, it includes
the second tuning parameter $\beta_{i}^{}$\,;
moreover, the algorithm (\ref{EGA}) uses initial value $\theta_0$
in a different manner;
\item along with (\ref{alg1}), our algorithm contains
the additional step of averaging of the updates (\ref{alg2}).
\end{itemize}

Papers \cite{HeKiWa99,KiWa01} study convergence properties of the EG
method (\ref{EGA}) in a deterministic setting. Namely, they show
that, under some assumptions, the difference $A_t(\theta_t) -
\min_{\theta\in\Theta_{M,1}}A_t(\theta)$ is bounded by a constant
depending on $M$ and $t$. If this constant is small enough, these
results show that the EG method provides good numerical minimizers
of the empirical risk $A_t$. However, they do not apply to the
expected risk. In particular, they do not imply that the expected
risk $\mathbb{E} A(\theta_t)$ is close to the minimal possible value
$\min_{\theta\in\Theta_{M,1}}A(\theta)$, which, as we prove, is true
for the algorithm with averaging proposed here.
%Simulations show
%that averaging of the updates is essential: the generalization error
%of the EG method without averaging is typically much larger than
%that of our method \cite{NIPS05}.

Finally, we point out that the algorithm (\ref{MDA}) may be deduced
from the ideas mentioned in Subsection \ref{sub:heur} and which are
studied in the literature on proximal methods within the field of
convex optimization (see, e.g., \cite{Ki97,BeTe03}  and the
references therein). Namely, under rather general conditions, the
variable $\theta_i$ from (\ref{MDA}) is the solution to the
minimization problem
\begin{equation}\label{hheur}
\theta_i=\argmin_{\theta\in \Theta}\left( \theta^T \gamma_i u_i(\theta_{i-1})+
B(\theta,\theta_{i-1})\right),
\end{equation}
where the penalty $B(\theta,\theta_{i-1}) = V(\theta) -
V(\theta_{i-1})- (\theta-\theta_{i-1})^T\nabla V(\theta_{i-1})$
represents the Bregman divergence between $\theta$ and
$\theta_{i-1}$ related to the strongly convex function $V$.

\section{Proofs}\label{sec:rate}

In this section, we provide technical details leading to the result
of Theorem \ref{thm}. They will be given in a more general setting
than that of Theorem \ref{thm} (cf Theorem 1 of \cite{Nest05-1}). Namely, we will consider an
arbitrary proxy function $V$ and use the notations and assumptions
of Subsection \ref{sub:algo}. Propositions \ref{prop:bound1} and
\ref{prop:bound2} below are valid for an arbitrary closed convex
set $\Theta$ in $E$, and for the estimate sequences $(\theta_i)$
and $(\widehat{\theta}_i)$ defined by the algorithm
(\ref{alg1})--(\ref{alg2}). The argument up to the relation
(\ref{ai2}) in the proof of Theorem \ref{thm} is valid under the
assumption that $\Theta$ is a convex compact set in $E$.

\medskip

Introduce the notation
\begin{eqnarray*}
\forall\, \theta\in \Theta, \quad    \nabla A (\theta) &=& \mathbb{E}\, u_i(\theta), \\
 \xi_i(\theta) &=& u_i(\theta) - \nabla A (\theta)~,
\end{eqnarray*}
where the random functions $u_i(\theta)$ are defined in (\ref{1}).
Note that the mapping $\theta\mapsto \mathbb{E}\, u_i(\theta)$
belongs to the sub-differential of $A$ (which explains the notation
$\nabla A$). This fact and the inequality
$\mathbb{E}\,\|u_i(\theta)\|_\infty^2\le L_\varphi^2(\lambda)$, valid
for all $\theta\in \Theta$, are the only properties of $u_i$ that
will be used in the proofs, other specific features of definition
(\ref{1}) being of no importance.

\begin{proposition}\label{prop:bound1}
For any $\theta\in\Theta$ and any integer $t\ge 1$ the following
inequality holds
\begin{eqnarray*}
   && \sum_{i=1}^{t} \gamma_i (\theta_{i-1}-\theta)^T
    \nabla A(\theta_{i-1}) \\
   &&\le \beta_t V(\theta) - \beta_0 V(\theta_*)
 -\sum_{i=1}^{t}\gamma_i (\theta_{i-1}-\theta)^T \xi_i(\theta_{i-1})
 +\sum_{i=1}^{t}
    \frac{\gamma_i^2}{2\alpha\beta_{i-1}} \|u_i(\theta_{i-1})\|_\infty^2\,.
\end{eqnarray*}
\end{proposition}

\medskip

\begin{proof}
By continuous differentiability of $W_{\beta_{t-1}}$, we have
\[
W_{\beta_{i-1}}(\zeta_i) = W_{\beta_{i-1}}(\zeta_{i-1}) + \int_0^1
(\zeta_i-\zeta_{i-1})^T \nabla W_{\beta_{i-1}}(\tau\zeta_i
+(1-\tau)\zeta_{i-1})\, {\rm d}\tau \, .
\]
Put $v_i = u_i(\theta_{i-1})$. Then $\zeta_i-\zeta_{i-1}=\gamma_i v_i$,
and by the Assumption (L)
\begin{eqnarray*}
W_{\beta_{i-1}}(\zeta_i) &=& W_{\beta_{i-1}}(\zeta_{i-1})+ \gamma_iv_i^T\nabla
W_{\beta_{i-1}}(\zeta_{i-1})
\\
 &&+ \ \gamma_i \int_0^1 v_i^T\Big[
\nabla W_{\beta_{i-1}}(\tau\zeta_i +(1-\tau)\zeta_{i-1})-\nabla
W_{\beta_{i-1}}(\zeta_{i-1})\Big] {\rm d}\tau
\\
&\le& W_{\beta_{i-1}}(\zeta_{i-1})+ \gamma_iv_i^T\nabla
W_{\beta_{i-1}}(\zeta_{i-1})
\\
 &&+ \ \gamma_i \|v_i
\|_\infty^{} \int_0^1 \| \nabla W_{\beta_{i-1}}(\tau\zeta_i
+(1-\tau)\zeta_{i-1})-\nabla W_{\beta_{i-1}}(\zeta_{i-1})\|_1^{} {\rm d}\tau
\\
&\le& W_{\beta_{i-1}}(\zeta_{i-1})+ \gamma_iv_i^T\nabla
W_{\beta_{i-1}}(\zeta_{i-1}) + \frac{\gamma_i^2 \|v_i
\|_\infty^2}{2\alpha\beta_{i-1}}\;.
\end{eqnarray*}
Using the last inequality, the fact that $(\beta_i)_{i \ge 1}$ is a
non-decreasing sequence and that, for $z$ fixed, $\beta \mapsto
W_{\beta}(z)$ is a non-increasing function, we get
\[
W_{\beta_i}(\zeta_i) \le W_{\beta_{i-1}}(\zeta_i) \le
W_{\beta_{i-1}}(\zeta_{i-1})
 -\gamma_i \theta^T_{i-1}v_i
  + \frac{\gamma_i^2 \|v_i \|_\infty^2}{2\alpha\beta_{i-1}}\,.
\]
When summing up, we obtain
\[
\sum_{i=1}^{t} \gamma_i \theta_{i-1}^T v_i \le W_{\beta_0}(\zeta_0)
-W_{\beta_t}(\zeta_t)
  + \sum_{i=1}^{t}
    \frac{\gamma_i^2 \|v_i \|_\infty^2}{2\alpha\beta_{i-1}}\,.
\]
Using the representation $\zeta_t=\sum_{i=1}^{t}\gamma_i v_i$ we get
that, for any $\theta\in\Theta$,
\[
\sum_{i=1}^{t} \gamma_i (\theta_{i-1}-\theta)^T v_i \le W_{\beta_0}(\zeta_0) -
W_{\beta_t}(\zeta_t) -\zeta_t^T\theta
  + \sum_{i=1}^{t}
    \frac{\gamma_i^2 \|v_i \|_\infty^2}{2\alpha\beta_{i-1}} \,.
\]
Finally, since $v_i = \nabla A(\theta_{i-1}) + \xi_i(\theta_{i-1})$
we find
\begin{eqnarray*}
&&\sum_{i=1}^{t} \gamma_i (\theta_{i-1}-\theta)^T
    \nabla A(\theta_{i-1})
\\
&&\le W_{\beta_0}(\zeta_0)-W_{\beta_t}(\zeta_t) -\zeta_t^T\theta
    -\sum_{i=1}^{t} \gamma_i (\theta_{i-1}-\theta)^T \xi_i(\theta_{i-1})
  + \sum_{i=1}^{t}
    \frac{\gamma_i^2 \|v_i \|_\infty^2}{2\alpha\beta_{i-1}} \, .
\end{eqnarray*}
Thus, the desired inequality follows from the fact that
\[
W_{\beta_0}(\zeta_0)=W_{\beta_0}(0) = \beta_0
\sup_{\theta\in\Theta}\{-V(\theta)\} = -\beta_0 V(\theta_*)
\]
and $\beta V(\theta) \ge -W_\beta(\zeta) -\zeta^T\theta$,  for all
$\zeta\in \RR^M$.
\end{proof}

\medskip

Now we derive the main result of this section.

\medskip

\begin{proposition}\label{prop:bound2}
For any integer $t\ge 1$, the following inequality holds true:
%\begin{eqnarray}%\nonumber
\begin{equation}
\mathbb{E}\, A(\widehat{\theta}_t) \le
\inf_{\theta\in\Theta}\left[ A(\theta) + \displaystyle
\frac{\beta_t V(\theta)- \beta_0 V(\theta_*)}{\sum_{i=1}^t
\gamma_i} \right]
%\\
 + \
%&&
L^2_{\varphi}(\lambda)\left(\sum_{i=1}^t \gamma_i\right)^{-1}
    \sum_{i=1}^t\frac{\gamma_i^2}{2\alpha\beta_{i-1}}
    \, .\label{ora}
\end{equation}
%\end{eqnarray}
Hence, the expected accuracy of the estimate $\widehat{\theta}_t$
with respect to the $\varphi$-risk satisfies the following upper
bound:
\begin{equation}\label{main}
\mathbb{E}\, A(\widehat{\theta}_t) - \min_{\theta\in\Theta}
A(\theta) \le \displaystyle \frac{1}{\sum_{i=1}^t \gamma_i}
\left(\beta_t V(\theta^*_A)- \beta_0 V(\theta_*)\ + \
L^2_{\varphi}(\lambda)
    \sum_{i=1}^t\frac{\gamma_i^2}{2\alpha\beta_{i-1}}
            \right)
    \, ,
\end{equation}
where $\theta^*_A \in\Argmin_{\theta\in\Theta} A(\theta)$.
\end{proposition}
\medskip

\begin{proof}
For any $\theta\in\Theta$, by convexity of $A$, we get
\begin{eqnarray}\nonumber
\mathbb{E}\, A(\widehat{\theta}_t) - A(\theta)
 &\le& \frac{\sum_{i=1}^t \gamma_i \left( \mathbb{E}\, A(\theta_{i-1}) -
A(\theta)\right)}{\sum_{i=1}^t \gamma_i} \nonumber
\\
 &\le& \frac{\sum_{i=1}^t \gamma_i \mathbb{E}\,
 [\left( \theta_{i-1} - \theta
                                                \right)^{T}
 \nabla A(\theta_{i-1})]}{\sum_{i=1}^t \gamma_i}\,.
 \label{3}
\end{eqnarray}
Conditioning on $\theta_{i-1}$ and then using both the definition
of $\xi_i(\theta_{i-1})$ and the independence between
$\theta_{i-1}$ and $(X_i,Y_i)$, we obtain
\[
\mathbb{E}\, \xi_i(\theta_{i-1}) = 0~.
\]
We now combine (\ref{3}) and the inequality of Proposition
\ref{prop:bound1} where we make use of the bound
$\mathbb{E}\,\|u_i(\theta)\|_\infty^2\le L_\varphi^2(\lambda)$. This
leads to (\ref{ora}). Inequality (\ref{main}) is straightforward in
view of (\ref{ora}).
\end{proof}

\medskip

\begin{remark}
Note, that simultaneous change of scale in the definition of the
sequences $(\beta_i)$ and $(\gamma_i)$  (i.e. multiplying them by
the same positive constant factor) does not affect the upper bounds
in the previous propositions, though it might affect the estimate
sequences $(\theta_i)$ and $(\widehat{\theta}_i)$ of the algorithm
(\ref{alg1})--(\ref{alg2}).
\end{remark}

\medskip

\noindent {\bf Proof of Theorem \ref{thm}.} We have
$V(\theta_*)=0$, and
\[
V(\theta^*_A)\le \max_{\theta\in\Theta} V(\theta)\triangleq V^*.
\]
Using (\ref{main}) with the choice $\gamma_i\equiv 1$
and $\beta_i=\beta_0\sqrt{i+1}$ for
$\beta_0>0$, $i\ge 1$,
we get
\begin{equation}\label{ai2AVN}
\mathbb{E}\, A(\widehat{\theta}_t) - A(\theta^*_A)
 \le
%{1\over \sqrt{t}}
 {\sqrt{t+1}\over t}
 \left(\beta_0  V^*+{L^2_{\varphi}(\lambda)
 \over \alpha\beta_0}\right).
\end{equation}
Optimizing this bound in $\beta_0$ leads to the choice:
\[
\beta_0=\frac{L_{\varphi}(\lambda)}{\sqrt{\alpha V^*}}\ ,
\]
which gives the bound:
\begin{equation}
\label{ai2} \mathbb{E}\, A(\widehat{\theta}_t) - A(\theta^*_A)
 \le
 \frac{2L_{\varphi}(\lambda)}{t}\sqrt{\frac{V^*}{\alpha}\,(t+1)}\,.
% 2L_{\varphi}(\lambda)\,\sqrt{\frac{V^*}{\alpha t}}\,.
\end{equation}

We now recall that $\Theta= \Theta_{M,\,\lambda}$ and that, for
the proxy function $V$ defined in (\ref{2}), we have $\alpha=
\lambda^{-1}$. Furthermore, this proxy function attains its
maximum at each vertex of the $\lambda$-simplex
$\Theta_{M,\,\lambda}$ and satisfies
\begin{equation}\label{d1}
V^*=\max_{\theta\in\Theta_{M,\,\lambda}} V(\theta) = \lambda \ln M \, .
\end{equation}
Therefore, the optimal value $\beta_0$ equals $
L_{\varphi}(\lambda) (\ln M)^{-1/2}$. This gives the accuracy
bound as in the statement of Theorem \ref{thm}.\
\rule{0.5em}{0.5em}

\section{Extension}\label{sec:gener}

Theorem \ref{thm} can easily be extended to a more general
framework. Inspection of the proof indicates that it does not use a
specific form of the loss function or of the proxy function. The
required properties of these functions are summarized at the
beginning of Section \ref{sec:rate}. We now state a more general
result. First, introduce some notation.

Consider a random variable $Z$ which takes its values in a set
$\mathcal{Z}$. The decision set $\Theta$ is supposed to be a convex
and closed set in $\mathbb{ R}^M$, and a loss function
$Q:\Theta\times\mathcal{Z}\to\mathbb{R}_+$ such that the random
function\, $Q(\cdot\,,Z):\Theta\to\mathbb{R}_+$\, is convex almost
surely. Define the convex risk function $A:\Theta\to\mathbb{R}_+$ as
follows:
\[
A (\theta) = \mathbb{E}\, Q(\theta,Z)\,.
\]
The training sample is given in the form of an i.i.d. sequence
$(Z_1,\dots,Z_{t-1})$, where each $Z_i$ has the same distribution
as $Z$. Our aim now consists in \textit{criterial minimization} of
$A$ over $\Theta$ (see, e.g., \cite{PolTsyp84}), which means that
we characterize the accuracy of the estimate
$\widehat{\theta}_t=\widehat{\theta}_t(Z_{1},\dots,Z_{t-1})\in\Theta$
minimizing $A$, by the difference:
\[
\mathbb{E}\, A (\widehat{\theta}_t) -\min_{\theta\in\Theta} A
(\theta)
\]
(we assume that $\displaystyle\min_{\theta\in\Theta}A(\theta)$ is
attainable). We denote by
\begin{equation}\label{Genu}
u_i(\theta) = \nabla_{\theta}^{}Q(\theta,Z_i )\,, \quad i=1,2,\dots,
\end{equation}
the stochastic sub-gradients which are measurable functions
defined on $\Theta\times\mathcal{Z}$ such that, for any
$\theta\in\Theta$, their expectation $\mathbb{E}\,u_i(\theta)$
belongs to the sub-differential of the function $A(\theta)$.

\medskip

\begin{theorem}\label{th2}
Let $\Theta$ be a convex closed set in $\mathbb{R}^M$, and $Q$ be a
loss function which meets the conditions mentioned above. Moreover,
assume that:
\begin{equation}
\label{eq:supth2} \sup_{\theta\in\Theta}\mathbb{E}\,
\|\nabla_{\theta}^{}Q(\theta,Z)\|_{\infty}^{2}\leq L_{\Theta,Q}^2~,
\end{equation}
where $L_{\Theta,Q}$ is a finite constant. Let $V$ be a proxy
function on $\Theta$ satisfying Assumption (L) with a parameter
$\alpha>0$, and assume that there exists $\theta^*_A
\in\Argmin_{\theta\in\Theta} A(\theta)$. Then, for any integer $t\ge
1$, the estimate $\widehat{\theta}_t$, defined in Subsection
\ref{sub:algo} with stochastic sub-gradients (\ref{Genu}) and with
sequences $(\gamma_i)_{i\ge 1}$ and $(\beta_i)_{i\ge 1}$ from
(\ref{gb1}) with arbitrary $\beta_0 >0$, satisfies the following
inequality:
\begin{eqnarray}
\label{eq:th3} \mathbb{E}\, A(\widehat{\theta}_t) -
\min_{\theta\in\Theta}A(\theta)
 \le
% \frac{1}{\sqrt{t}}
 \left(\beta_0  V(\theta^*_A)+{L_{\Theta,Q}^2
 \over \alpha\beta_0}\right)
 {\sqrt{t+1}\over t}
 \,.
\end{eqnarray}
Furthermore, if $\overline{V}$ is a constant such that
$V(\theta^*_A)\le \overline{V}$ and we set $\beta_0 =
L_{\Theta,Q}\,(\alpha\, \overline{V}\,)^{-1/2}$, then
\begin{eqnarray}
\label{eq:th31} \mathbb{E}\, A(\widehat{\theta}_t) -
\min_{\theta\in\Theta}A(\theta)
 \le
 2\, L_{\Theta,Q}
 \left(\alpha^{-1}\overline{V}\right)^{1/2} {\sqrt{t+1}\over t}
% \sqrt{\frac{\overline{V}}{\alpha t} }
\, .
\end{eqnarray}
In particular, if\, $\Theta$ is a convex compact
%subset of $\mathbb{R}^M$,
 set, we can take $\overline{V} = \max_{\theta\in\Theta}
V(\theta)$.
%for $\beta_0 =  L_{\Theta,Q}(\alpha \, \displaystyle
%\max_{\theta\in\Theta} V(\theta))^{-1/2}$ we have
%\begin{eqnarray}
%\label{eq:th2} \mathbb{E}\, A(\widehat{\theta}_t) -
%\min_{\theta\in\Theta}A(\theta)
% \le
% 2\, L_{\Theta,Q}
% \left(\alpha^{-1}\max_{\theta\in\Theta} V(\theta)\right)^{1/2}
% {\sqrt{t+1}\over t} \, .
% \end{eqnarray}
\end{theorem}

This theorem follows from the proofs of Section \ref{sec:rate} (cf.
(\ref{main}),  (\ref{ai2AVN}), and (\ref{ai2})), where
$L_{\Theta,Q}$ should replace the constant $L_{\varphi}(\lambda)$.
It generalizes Theorem \ref{thm} and encompasses different
statistical models, including the one described in Section
\ref{sec:ProblemResult}, where $Z$ plays the role of the pair of
variables $(X,Y)$, $\Theta=\Theta_{M,\,\lambda}$, and
$Q(\theta,Z)=\varphi(Y\theta^TH(X))$. In the same way, Theorem
\ref{th2} is also applicable to the standard regression model with
squared loss $Q(\theta,Z)=(Y-\theta^TH(X))^2$, in which case a
similar result has been proved for another method in \cite{JuNe00}.

\begin{remark}({\sc dependent data}). Inspection of the proofs shows
that Theorem \ref{th2} can be easily extended to the case of
dependent data $Z_i$. In fact, instead of assuming that $Z_i$ are
i.i.d., it suffices to assume that they form a stationary sequence,
where each $Z_i$ has the same distribution as $Z$. Then Theorem
\ref{th2} remains valid if we additionally assume that the
conditional expectation $\mathbb{E}\,
(\xi_i(\theta_{i-1})|\theta_{i-1}) = 0$ a.s.
%This assumption is verified, for example, for stationary
%autoregression model with squared loss.
\end{remark}

%Finally note, that constant $L_{\Theta,Q}^{}$ from Theorem~2 may be raised,
%Maiming at its simplier calculation, by the supremum in (\ref{eq:supth2})
%to be taken over a wider subset of $\Theta$, for instance, in case of
%$\Theta\subset\{\theta\,:\, \|\theta\|_1^{}\leq\lambda\}$ for some
%$\lambda>0$; in a similar manner, one may also deal with calculation
%of constant $\overline V\geq \displaystyle\max_{\theta\in\Theta}V(\theta)$.
%\bigskip

\section{Discussion}\label{sec:conclu}

To conclude, we discuss further the choices of the proxy function
$V$, of the parametric set $\Theta$ and of the sequences
$(\beta_i)_{i\geq 1}$, $(\gamma_i)_{i\geq 1}$.

\subsection{Choice of the proxy function $V$}\label{sub:proxy-fu}

The choice of the entropic proxy function defined in (\ref{2}) is
not the only possible. A key condition on $V$ is the strong
convexity with respect to the norm $\|\cdot\|_1^{}$, which
guarantees that Assumption (L) holds true.
Therefore, one may also consider other
proxy functions satisfying this condition, such as, for example:
\begin{equation}\label{Vst2}
    \forall\, \theta\in\RR^M\,, \quad
V(\theta) = \frac{1}{2\lambda^2} \|\theta\|_p^2
   = \frac{1}{2\lambda^2} \left( \sum_{j=1}^M \left(\theta^{(j)} \right)^p
                          \right)^{2/p} ,
\end{equation}
where $p=1+1/\ln M$\, (see \cite{BeNe99}). In contrast with the
function (\ref{2}), the proxy function (\ref{Vst2}) can be used when
$\Theta$ is any convex and closed set in $\RR^M$.

For the simplex $\Theta_{M,\,\lambda}$, one may consider functions
of the form
\begin{equation}\label{Vst3}
    \forall\, \theta\in\Theta_{M,\,\lambda}\,, \quad
V(\theta) = C_0 + C_1\sum_{j=1}^M \left(\theta^{(j)}\right)^{s+1} ,
\quad s=\frac{1}{\ln M}\,,
\end{equation}
where the constants $C_0= - \lambda^{2}/({\rm e} s(s+1))$, $C_1=
\lambda^{1-s}/(s(s+1))$ are adjusted in order to have
$\min_{\theta\in\Theta_{M,\,\lambda}}V(\theta)=0$. It is easy to see
that the proxy function defined in (\ref{Vst3}) is $\alpha$-strongly
convex in the norm $\|\cdot\|_1^{}$. When $\lambda=1$, this proxy
function equals to a particular case of $f$-divergence of Csisz\'ar
(see the definition in \cite{Vai86}) between the uniform
distribution on the set $\{1,\dots,M\}$ and the distribution on the
same set defined by probabilities $\theta^{(j)}$. Recall that for
$\lambda=1$ the proxy function defined in (\ref{2}) equals to the
Kullback divergence between these distributions. Presumably,
other proxy functions can be based on some properly chosen
$f$-divergences of Csisz\'ar.

On the other hand, if a proxy function $V$ is such that the gradient
of its $\beta$-conjugate $\nabla W_{\beta}$ cannot be explicitly
written, the numerical implementation of our algorithm might become
time-consuming.

>From the upper bound (\ref{eq:th31}) we can see that an important
characteristic of $V$ is the ratio $\overline{V}/\alpha$ (or
$(\max_{\theta\in\Theta}V(\theta)/\alpha)$ if the set $\Theta$ is
bounded) and thus one can look for optimal proxy functions
minimizing this ratio. We conjecture that for
$\Theta=\Theta_{M,\,\lambda}$ such an optimal proxy function is the
entropy type function given in (\ref{2}); however, we do not have a
rigorous proof of this fact. For the latter, we have
\[
\frac{1}{\alpha} \max_{\theta\in\Theta_{M,\,\lambda}} V(\theta)
= \lambda^2 \ln M \,.
\]
For other proxy functions, this ratio is of the same order. For
instance, it is proved in \cite{BeNe99} (Lemma 6.1) that the proxy
function defined in (\ref{Vst2}) satisfies
\[
\frac{1}{\alpha} \max_{\theta\in\Theta_{M,\,\lambda}} V(\theta)
= O(1) \lambda^2 \ln M \,.
\]
This relation is true for the proxy function defined in (\ref{Vst3})
as well.

%Finally, note that such widely used penalty functions as $\| \cdot
%\|_1$ and $\| \cdot \|_2^{2}$ are not proxy functions in the sense
%of Definition \ref{proxyf}: they do not satisfy Assumption~(L).
Finally, note that a widely used penalty function as $\| \cdot \|_1$
is not a proxy function in the sense of Definition \ref{proxyf} as
it is not strongly convex with respect to $\|\cdot\|_1$. Another
frequently used penalty function %--- the function
$V(\theta)=\|\theta\|^2_2$ %---
is strongly convex. However, it can be
easily verified that its ``performance ratio" is extremely bad for large $M$:
this function $V$ satisfies
\[
\frac{1}{\alpha} \max_{\theta\in\Theta_{M,\,\lambda}}V(\theta)
=\frac{1}{2}\lambda^2M\,.
\]

\subsection{Other parametric sets $\Theta$}\label{sub:set}

%Theorem \ref{th2} holds for any convex closed $\Theta$ contained in
%$\RR^M$. However, for general sets, the gradient $\nabla W_{\beta}$
%might be impossible to compute explicitly and then the algorithm
%will be hard to implement. Hence, it is important  to consider only
%relatively simple sets $\Theta$ such as: (i) the $\lambda$-simplex
%$\Theta_{M,\,\lambda}$, (ii) the full-dimensional $\lambda$-simplex
%$\left\{\theta\in\mathbb{R}_{+}^M: \|\theta\|_1^{}\leq \lambda
%%\sum_{i=1}^M \theta^{(i)}\leq\lambda\,, \quad \theta^{(j)}\geq 0,\, j=1,\dots,M
%\right\}$ and (iii) the symmetrized version of the latter, that is
%the hyper-octahedron
%$\{\theta\in\mathbb{R}^M:\|\theta\|_1^{}\leq\lambda\}$. This
%simplicity requirement on the set $\Theta$ can be removed if we are
%able to compute efficiently the anti-gradient $-\nabla
%W_{\beta_i}(\zeta_i)$ along the iterations (by solving the
%maximization problem stated in the assertion 2(i) of Proposition~\ref{prop:conj}).

Theorem \ref{th2} holds for any convex closed $\Theta$ contained in
$\RR^M$. However, for general sets, the gradient $\nabla W_{\beta}$
cannot be computed explicitly and the computation effort of
implementing an iteration of the algorithm can become prohibitive.
Hence, it is important to consider only the sets $\Theta$ for which
the solution $\theta_*(z)=-\nabla W_{\beta}(z)$ of the optimization
problem (\ref{fench}) can be easily computed. Some examples of such
``simple" sets are: (i) the $\lambda$-simplex
$\Theta_{M,\,\lambda}$, (ii) the full-dimensional $\lambda$-simplex
$\left\{\theta\in\mathbb{R}_{+}^M: \|\theta\|_1^{}\leq \lambda
%\sum_{i=1}^M \theta^{(i)}\leq\lambda\,, \quad \theta^{(j)}\geq 0,\, j=1,\dots,M
\right\}$ and (iii) the symmetrized version of the latter, that is
the hyper-octahedron
$\{\theta\in\mathbb{R}^M:\|\theta\|_1^{}\leq\lambda\}$.

\subsection{Choice of the step size and temperature parameters}
\label{sub:seque}

The constant factor in the bound of Theorem \ref{thm} can be only
slightly improved. The sequences $(\beta_i)$ and $(\gamma_i)$ as
described in (\ref{gb1}) are close to optimal ones in a sense of the
upper bound
 (\ref{main}). Indeed, if we further bound $V(\theta_A^*)$ by
 $V^*=\max_{\theta\in\Theta}V(\theta)$ in (\ref{main})
%$V(\theta_A^*)$ for $V^*=\max_{\theta\in\Theta}V(\theta)$
and minimize in $(\gamma_i)$ and $(\beta_i)$\, under the
monotonicity condition
 $\beta_i\geq\beta_{i-1}$, we get that the minimum is obtained for  sequences $(\gamma_i)$ and $(\beta_i)$,
which are independent of $i$ and such that $\beta_i/\gamma_i \equiv
L_\varphi(\lambda)\sqrt{t/2\alpha V^*}$. We can take, for instance,
\begin{equation}\label{Vst4}
\gamma_i\equiv\frac{1}{\sqrt{t}}\,,\quad \beta_i\equiv \frac{
L_\varphi(\lambda)}{\sqrt{2\alpha V^*}}\,,
\end{equation}
which leads to a better constant than the bound (\ref{eq:th1}) in
Theorem 1
\begin{equation}\label{optmult}
  \mathbb{E}\, A(\widehat{\theta}_t) - \min_{\theta\in\Theta_{M,\,\lambda}} A(\theta)
  \leq \lambda\, L_\varphi(\lambda) \sqrt{\frac{2\ln M}{t}}\,.
\end{equation}
Thus, we can improve the constant factor in the upper bound
 from $2$ to $\sqrt{2}$. However, in order to make this improvement, one needs to
know the sample size $t$ in advance, and this is not compatible with
the online framework.

%\newpage

\section*{Appendix}\label{appendix}

In this appendix, we propose two different proofs of the fact that
Assumption (L) holds for the $\beta$-conjugate $W_{\beta}$ of the
entropy type function $V$ given by (\ref{2}). First, we give a
straightforward argument using the equations from (\ref{gibbs}). The
second proof is based on a generic argument which exploits the
convexity properties of function $V$ rather than its particular
expression.

\subsection*{A. Direct proof}
Evidently, the function $W_{\beta}$ in (\ref{Wbeta}) is twice
continuously differentiable on $E^*=\ell_\infty^M$.
Set
%Let us denote
%$L$ the underlying Lipschitz constant. Note that
\begin{equation*}
  L = \sup_{z_1,z_2\in E^*,\,z_1\neq z_2}
        \frac{\|\nabla W(z_1)-\nabla W(z_2)\|_1}{\|z_1-z_2\|_\infty}
 \le \sup_{\|x\|_\infty=1, \, \|y\|_\infty\le 1}\;
         \sup_{z\in E^*}x^T\nabla^2 W(z)y
\end{equation*}
where the second derivative matrix $\nabla^2 W(z)$ has the entries
\[
\frac{\partial^2 W(z)}{\partial z_i\partial z_j}=
\frac{\lambda}{\beta} \left( {e^{-z_i/\beta}\delta_{ij}\over
\sum_ke^{-z_k/\beta}}-{e^{-z_j/\beta}e^{-z_i/\beta}\over
(\sum_ke^{-z_k/\beta})^2} \right).
\]
Here $\delta_{ij}$ stands for the Kronecker symbol.
Denote $a_i={e^{-z_i/\beta}/ \sum_ke^{-z_k/\beta}}$ which are
evidently positive with $\sum_ia_i=1$. Now,
\begin{eqnarray} \nonumber
\frac{\beta}{\lambda}\, x^T\nabla^2 W(z)y&=&\sum_{i}x_iy_i
a_i-\sum_{i}a_ix_i\sum_ja_jy_j=
\sum_{i}x_i a_i[y_i-\sum_ja_jy_j]\\
&=&\sum_{i}x_i a_i\sum_{j\neq i}a_j(y_i-y_j)\le \sum_{i}\sum_{j\neq
i}a_ia_j|y_i-y_j|\,. \label{yiyj}
\end{eqnarray}
Finally,
%We are done as
the latter sum is bounded by $1$ for any $|y_i|\le 1$ and $a_i\ge
0$, $\sum_ia_i=1$.
%\par
To see this, note that the maximum of the convex (in $y\in R^M$)
function of the right hand side (\ref{yiyj}) on the convex set
$\{y\in R^M :\|y\|_\infty\le 1\}$ is always attained at the extreme
points of the set, which are the vertices of the hypercube $\{y\in
R^M: y_i=\pm 1,\;i=1,\dots,M\}$. Denote this extreme point
$y^*=(y_1^{},\dots,y_M^{})^T$. Let us split the index set
$\{1,\dots,M\}$ into $I_+=\{i:\,y^*_i=1\}$ and
$I_-=\{i:\,y^*_i=-1\}$. Then the maximal value of the sum can be
decomposed, and we get
\[
\frac{\beta}{\lambda}\, L\le 2\sum_{i\in I_+, j\in I_-}a_ia_j
+2\sum_{j\in I_+, i\in I_-}a_ia_j= 4\sum_{i\in I_+} a_i\sum_{j\in
I_-}a_j=4a_+^{}(1-a_+^{})\le 1
\]
where $a_+^{}=\displaystyle\sum_{i=I_+} a_i$. Hence,
$\alpha=1/\lambda$ which is independent of
%both $M$ and
$\beta$.

\subsection*{B. General argument}

Consider first a more general setting, where $E$ is the space
$\mathbb{R}^M$ equipped with some norm $\|\cdot\|$ (primal space),
and $E^*$ is $\mathbb{R}^M$ equipped with the corresponding dual
norm $\|\cdot\|_*$ (dual space). Let now $V:\Theta\to\mathbb{R}$ be
$\alpha$-strongly convex function with respect to primal norm
$\|\cdot\|$ on a convex closed set $\Theta\subset\mathbb{R}^M$. %We
%assume that $\alpha$ is positive and does not depend on $M$
%(strictly speaking, we have here a family of functions $V$
%parameterized by integer $M\ge 2$).
It is easy to show that
inequality (\ref{sconv}) (holding for all $x,y\in\Theta$ and for any
$s\in[0,1]$) implies
\begin{equation}\label{x^{*}}
V(x)\ge V(x^*)+{\alpha\over 2}\|x-x^*\|^2,\quad\forall\,x\in\Theta,
\end{equation}
 where
$x^*=\argmin_{x\in \Theta} V(x).$ Indeed, the existence and
uniqueness of the minimizer $x^*$ is evident. Now for any $x\in
\Theta$ and $s\in (0,1)$,
\begin{eqnarray*}
  s V(x)+(1-s)V(x^*) &\ge& V(sx+(1-s)x^*)+\frac{\alpha}{2} s(1-s) \|x-x^*\|^2
 \\
   &\ge& V(x^*)+\frac{\alpha}{2} s(1-s) \|x-x^*\|^2
\end{eqnarray*}
and we get (\ref{x^{*}}) by subtracting $V(x^*)$, dividing by $s$,
and then letting $s$ tend to $0$ (as $V$ is continuous on $\Theta$).
\medskip

We assume, furthermore, that the $\beta$-conjugate $W_\beta$ defined
by (\ref{fench}) is continuously differentiable on $E^*$ and the
assertion 2(ii) of Proposition \ref{prop:conj} holds true.
%!
Let us fix any points $z_1,z_2\in E^*$ and arbitrary  $s\in(0,1)$.
%$0<s<1$.
Denote
\[
x_1=-\nabla W_\beta(z_1),\quad x_2=-\nabla W_\beta(z_2),
\]
and $x_s=s x_1+(1-s)x_2$\,. Recall that, due to the assertion 2(ii)
of Proposition \ref{prop:conj}, we have:
$x_k=\argmin_{\theta\in\Theta}\{z_k^T\theta+\beta V(\theta)\}$,
$k=1,2$.

Now we are done since the function $z^Tx+\beta V(x)$ is
$(\alpha\beta)$-strongly convex for any fixed $z$, hence, by
(\ref{x^{*}}),
\begin{eqnarray*}
-z_1^Tx_s-\beta V(x_s)&\le& -z_1^Tx_1- \beta V(x_1)-{\alpha\beta\over 2}\|x_s-x_1\|^2,\\
-z_2^Tx_s-\beta V(x_s)&\le& -z_2^Tx_2- \beta
V(x_2)-{\alpha\beta\over 2}\|x_s-x_2\|^2,
\end{eqnarray*}
and summing up with the coefficients $s$ and $1-s$ we get by
definition of $x_s$:
\begin{eqnarray*}
s(1-s)(z_1-z_2)^T(x_2-x_1)&=&s\,z_1^T(x_1-x_s)+(1-s)\,z_2^T(x_2-x_s)\\
&\le&\beta\left( V(x_s)-s V(x_1)-(1-s) V(x_2)
           \phantom{\alpha\over 2} \right.\\ &&\phantom{\beta(aaaaaaa)}\left.
                -{\alpha\over 2}\,s(1-s)\|x_1-x_2\|^2\right)\\
&\le& -\alpha\beta s(1-s)\|x_1-x_2\|^2.
\end{eqnarray*}
Therefore,
\[
\alpha\beta  \|x_2-x_1\|^2\le (z_1-z_2)^T(x_1-x_2)\le
\|z_1-z_2\|_*\,\|x_1-x_2\|\,,
\]
and it implies, both for $x_1=x_2$ and for $x_1\neq x_2$,
\[
\|\nabla W_\beta(z_1)-\nabla W_\beta(z_2)\| \le
\frac{1}{\alpha\beta} \|z_1-z_2\|_*
\]
which implies the desired Lipschitz property for $\nabla W_\beta$ in
the assertion 2(i) of Proposition \ref{prop:conj}.
\medskip

Now we return to the particular case where $\|\cdot\|=\|\cdot\|_1$,
$\|\cdot\|_*=\|\cdot\|_\infty$ and $V$ is the entropy type proxy
function $V$ defined in (\ref{2}). We prove that $V$ is
$(1/\lambda)$-strongly convex with respect to the norm
$\|\cdot\|_1$\,, i.e. it satisfies (\ref{sconv}) with
$\alpha=1/\lambda$.

\emph{Proof of  (\ref{sconv}).} Observe that function $V$ defined in
(\ref{2}) is twice continuously differentiable at any point
$x=(x_1^{},\dots,x_M^{})^T$ inside the set $\Theta_{M,\lambda}$,
with
\[
\frac{\partial^2 V(x)}{\partial x_i^{}\partial x_j^{}}
=\frac{\delta_{ij}}{x_i^{}} \,,\quad i,j=1,\dots,M\,.
\]
Let us fix two arbitrary points $x,y$ inside the set
$\Theta_{M,\lambda}$. One may write, for some interior point
$\widetilde{x}\in\Theta_{M,\lambda}$,
\begin{eqnarray}\label{conv1}
  \langle\, \nabla V(x) - \nabla V(y) , x-y  \rangle\,  &=& (x-y)^T \nabla V(\widetilde{x})(x-y)
%  = \sum_{i=1}^M \frac{(x_i^{}-y_i^{})^2}{\widetilde{x}_i^{}}
   \\ \nonumber
   &=& \lambda^{-1} \sum_{i=1}^M \frac{\widetilde{x}_i^{}}{\lambda}\,
        \left(\frac{|x_i^{}-y_i^{}|}{\widetilde{x}_i^{}/\lambda}\right)^2
    \\ \label{Jensen}
   &\ge& \lambda^{-1} \left(\sum_{i=1}^M \frac{\widetilde{x}_i^{}}{\lambda}\,
        \frac{|x_i^{}-y_i^{}|}{\widetilde{x}_i^{}/\lambda}
                \right)^2
    \\ \label{conv2}
   &=& \lambda^{-1} \|x-y\|_1^2
\end{eqnarray}
where we used Jensen's inequality in (\ref{Jensen}) since all
$\widetilde{x}_i>0$ and
\[
\sum_{i=1}^M \frac{\widetilde{x}_i^{}}{\lambda} =1\,.
\]
By the standard argument (see, e.g., \cite{RoWe98}), for all
interior points $x,y$ of $\Theta_{M,\lambda}$ we get (\ref{sconv})
from (\ref{conv1})--(\ref{conv2}). Finally, by
%by Lemma \ref{auxconexity}, we get the result.
continuity of $V$ on $\Theta_{M,\lambda}$, (\ref{sconv}) extends to
all $x,y$ in $\Theta_{M,\lambda}$. \;{\footnotesize$\blacksquare$}

%\newpage

\end{document}